\documentclass[12pt]{amsart}
\input{xypic}
\usepackage{amssymb}

\newtheorem{theorem}{Theorem}[section]
\newtheorem{lemma}[theorem]{Lemma}

\newtheorem{proposition}[theorem]{Proposition}

\numberwithin{equation}{section}

\newcommand{\ra}{\rightarrow}

\def\umono{\ar@{_{(}->}[u]}
\def\uumono{\ar@{_{(}->}[uu]}

\def\lmono{\ar@{_{(}->}[l]}
\def\llmono{\ar@{_{(}->}[ll]}

\newcommand{\Z}{{\mathbb Z}}
\newcommand{\F}{{\mathbb F}}

\def\iso{\cong}

\title[Idempotency of CWL]{On the idempotency of some composite functors}

\author{Ram\'on J. Flores}
\thanks{The author is supported by MEC grant MTM2004-06686, MEC, Spain.}

\date{\today}

\begin{document}

\begin{abstract}

We present examples of localization functors whose composition
with certain cellularization functors is not idempotent, and
vice versa.

\end{abstract}

\maketitle




\section{Introduction}

The main goal of this note is to discuss some interesting
questions raised by E.~Dror-Farjoun, concerning the effect of
the successive application of localization and cellularization
functors (see definitions in next section). We deal with Conjectures 3 and 7 of \cite[pages 62--63]{Mislin
08}, which essentially ask for any localization functor $L$ and
cellularization $CW$, which are homotopy idempotent by
definition, whether the composite functors $LCW$ and $CWL$ are also homotopy idempotent.
Note that the functors $L$ and $CW$ are respectively coaugmented and
augmented, but in general neither $LCW$ nor $CWL$ admit
coaugmentation nor augmentation.

In light of these considerations
we use a previous description \cite[Corollary
5.8]{Flores07} of the $B\Z /2$-cellularization of the smallest
in the family of Suzuki simple groups,
$Sz(8)$, to provide counterexamples to the questions of Farjoun.
The following theorem summarizes our main results, which appear in
the text as Proposition \ref{CWL} and Proposition \ref{LCW}:

\begin{theorem}

If $CW$ represents cellularization with regard to a wedge $B\Z
/2\vee\bigvee_{\textrm{p odd prime}} M(\Z /p,2)$ and $L$ denotes an
$n$-Postnikov section for an appropriate choice of $n$, then the
functors $LCW$ and $CWL$ are not homotopy idempotent.

\end{theorem}

The paper is structured as follows: next section is devoted to recall some notions
of (co)Localization Theory, and also to a little
comprehensive introduction of the problem of the idempotency of the functors $LCW$ and $CWL$; Section 3 contains
our main results and their proofs.

\emph{Acknowledgments}. We warmly thank Carles Broto, Richard Foote and J\'{e}r\^{o}me Scherer for their careful reading
of the paper and their valuable comments.


\section{Preliminaries}

Recall that given a pointed map
$f:A\ra B$, a space $X$ is said to be $f$-local
if the following conditions hold:
the induced function
$\textrm{map}_*(B,X)\ra \textrm{map}_*(A,X)$ is a weak
equivalence; the localization associated to $f$ is a
coaugmented and
idempotent endofunctor $L$ of the category of spaces;
and every map $X\ra Y$ from a space $X$ to a
$f$-local space $Y$ factors through $LX$ in an unique way, up to
homotopy. Analogously, given spaces $A$ and $X$, $X$ is said
$A$-cellular if it is homotopy equivalent to an iterated homotopy
colimit of copies of $A$.  The $A$-cellularization of $X$,
denoted by $CW_A$, is
defined as an augmented idempotent endofunctor of the category of
spaces such that the augmentation $CW_AX\ra X$ is initial among
all maps $f:Y\ra X$ which induce a weak equivalence
$\textrm{map}_*(A,Y)\ra\textrm{map}_*(A,X)$. A comprehensive
introduction
to these functors can be found in the first two chapters of
\cite{Farjoun95}.

It is quite easy to find pairs of functors $CW$ and $L$ such that
$CWL$ or $LCW$ are idempotent. For example, if $L$ is a Postnikov
section $P[n]$, or more generally localization with regard the
constant map $A\ra {*}$ for a certain space $A$ (usually called
$A$-nullification), and $CW$ is the $n$-connected cover functor
(respectively $CW_A$), both composites are trivial. A similar
phenomenon occurs if we apply $CW_{M(\Z/p,2)}$, where $p$ is a
prime, and then take the Bousfield-Kan $q$-completion at another prime $q$.

It can also happen that the composite of $L$ and $CW$ behaves like the identity functor
(up to homotopy) over certain spaces, which immediately implies idempotency. An easy
example is given by taking $L$ as the $R$-completion, for
$R=\mathbb{Q}$ or $R=\Z /p$, of a simply-connected space, and $CW$
as the universal cover functor; then $LCW$ is the identity over $LX$, and 
$CWL$ is idempotent over $X$. So the problem is to find
functors $L$ and $CW$ such that their combined effect is not so
drastic as to kill the entire space, nor so weak that the composition
gives back the original space. As one can expect, for ``nice"
spaces homological localization functors \cite{Bousfield75}
combined with cellularization with regard to $p$-primary
meaningful spaces (as $B\Z /p$ or Moore spaces for $\Z /p$) have a
good behavior in this context, and idempotency can be expected at
least when working with nilpotent spaces. The same happens for some
well-known nullification functors
(\cite{Bousfield94},\cite{Bousfield97}).

So, when searching for counterexamples there are two different
but complementary paths: to seek for exotic
localizations/cellularizations, or to apply the aforementioned
functors to non-nilpotent spaces. The first strategy seems rather
unpromising, as precise descriptions of the effect of $L$, or
even more so of $CW$, are in general not available for functors that
are outside the context of $p$-primary homotopy; and even less
is known of
the structure of $LCWX$ or $CWLX$.
If the second path is taken, one quickly realizes that two
features of $L$ and $CW$ are needed: their effect on the target
space $X$ should be very precisely known, and both must change $X$
in a ``moderately opposite" way, so $CWLX$ or $LCWX$ should be
non-trivial but not very complicated spaces. This is the approach we take in next section, 
where our main results are proven.

\section{Idempotency and classifying spaces}

We will always consider cellularization with regard to the
classifying space $B\Z /2$, and we will denote this
functor simply by $CW$. In turn, if $M$ is the wedge of all Moore
spaces $M(\Z /q,2)$ for $q$ odd, $L$ will be nullification with
regard to a wedge of $M$ with a sphere $S^{n+1}$ whose dimension
will be determined in due course.

To keep this exposition brief we refer the reader to
\cite[16.4]{Gor80} for details on the 2-local structure of
$Sz(8)$.  As this group is our main
object of study, we henceforth denote it simply by $G$.

Before applying the functors, we need some technical remarks
concerning to the primary structure of $BG$. Bousfield-Kan
$q$-completion of a space $X$ (see \cite{Bousfield72} for a main
reference) is denoted by $X^{\wedge}_q$.

\begin{lemma}

The space $BG^{\wedge}_q$ is 2-connected for $q$ odd.

\label{2connected}
\end{lemma}

\begin{proof}

As $G$ is a simple group, the $q$-completion of $BG$ is
simply-connected, and according to \cite[VII.4.3]{Bousfield72},
$\pi_2BG^{\wedge}_q$ must be a finite $q$-group for every $q$. Moreover,
the Schur multiplier of $G$ is isomorphic to $\Z /2\times\Z /2$, so
$H_2(BG;\mathbb{F}_q)=0$ by universal coefficients. By Hurewicz's Theorem, we
have $\pi_2BG^{\wedge}_q=H_2(BG^{\wedge}_q;\mathbb{F}_q)$, which
is in turn isomorphic to $H_2(BG;\mathbb{F}_q)$ because $BG$ is
$q$-good \cite[VII.5]{Bousfield72}.
This completes the proof.

\end{proof}

Note that according to \cite[Theorem 6.2]{Rodriguez01}, the space
$BG^{\wedge}_q$ should be $M(\Z /q,2)$-cellular for every $q$ odd,
and then it is killed by $M(\Z /q,2)$-nullification, i.e.,
$\mathbf{P}_{M(\Z /q,2)}BG^{\wedge}_q$ is contractible.


Recall now from \cite{Gor80} that a Sylow 2-subgroup, $S$, of $G$
satisfies: $Z(S) = S' = \Omega_1(S) \iso (\Z /2)^3$ and
$S/Z(S) \iso (\Z /2)^3$.
Thus $Z(S) = Cl(S)$ is the minimal strongly closed subgroup that contains all
elements of order 2 in $S$.
The normalizer, $N$, of $S$ in $G$ acts transitively on the nonidentity
elements of both $Z(S)$ and $S/Z(S)$.
Following the
notation of \cite{Flores07}, let $\Gamma = N/Z(S) \iso (\Z /2)^3 \rtimes \Z /7$.
The following result
describes the low-dimensional homology of $\Gamma$:

\begin{lemma}

The group $\Gamma$ is 2-superperfect, i.e. $H_i(\Gamma,\Z /2)=0$ for
$i=1,2$. \label{superperfect}

\end{lemma}

\begin{proof}

The group $\Gamma$ is 2-perfect by \cite[page 51]{Flores07}, and
moreover it is not hard to check that $H_1(\Gamma;\Z)=\Z /7$.
Using now Gasch\"{u}tz's Theorem and the Fitting Lemma \cite{Gor80}, it can also be seen that the Schur
multiplier of $\Gamma$ is trivial. Then $H_2(\Gamma;\F_2)$ is
likewise trivial
by universal coefficients, and we are done.


\end{proof}

We can now give the precise definition of $L$. The group $\Gamma$ is
2-perfect and moreover has torsion in odd primes, so according to
\cite[1.1.4]{Levi95} $B\Gamma^{\wedge}_2$ has an infinite number
of nonzero homotopy groups. If $n$ is the smallest natural number
such that $\pi_nB\Gamma^{\wedge}_2$ is nontrivial, then our
localization $L$ will be nullification with regard to the wedge
$S^{n+1}\vee M$, where $M$ is the wedge of Moore spaces defined
above. Observe that according to the previous lemma and
\ref{2connected}, $n>2$.


We are ready to explore the behavior of $BG$ under the recursive
action of the functors $L$ and $CW$. As $n>1$, it is clear that
$BG$ is $L$-local, and therefore by \cite[Corollary
5.8]{Flores07}, $CWBG\simeq CWLBG$ fits in a fibre sequence:
$$
CWBG \longrightarrow BG \longrightarrow B\Gamma^{\wedge}_2\times\prod_{q\neq 2}
BG^{\wedge}_q
$$
where the product on the right is extended to odd
primes. Here the map on the right is given by the composition
$$
BG\longrightarrow BG^{\wedge}_2\longrightarrow BN^{\wedge}_2 \longrightarrow B\Gamma^{\wedge}_2,
$$
where the second map is an equivalence because the normalizer of
$S$ controls $G$-fusion in $S$ (see details in \cite[Example
5.2]{Flores07}).

Now let us describe $LCWBG$. Looping the base of the previous fibre
sequencee we obtain another one, which is principal:
\begin{equation} \label{looping}
\Omega (B\Gamma^{\wedge}_2)\times \Omega (\prod_{q\neq 2}
BG^{\wedge}_q) \longrightarrow CWBG \longrightarrow BG.
\end{equation}
The space $BG$ is clearly $S^{n+1}\vee M$-null, so by
\cite[Proposition 3.D.3]{Farjoun95} the fibration is preserved by
$L$. On the other hand, Lemma \ref{2connected} implies that
$\Omega (\prod_{q\neq 2} BG^{\wedge}_q)$ is simply connected, so
it is $M$-cellular by \cite[Theorem 6.2]{Rodriguez01} and its
$L$-localization is contractible. Then, by definition of $L$ we
have that $L\Omega (B\Gamma^{\wedge}_2)$ has the homotopy type of
an Eilenberg-Mac Lane space $K(H,n-1)$, where $H=\pi_n
B\Gamma^{\wedge}_2$. Thus $LCWBG(=LCWLBG)$ is homotopy equivalent to a
Postnikov piece $X$ whose only nonzero homotopy groups are
$\pi_1X=G$ and $\pi_{n-1}X=H$.

We are now ready to give a solution (in the negative) to the first
part of Problem~3 and Problem~7 in \cite[pages. 62-63]{Mislin 08}.

\begin{proposition}

The functor $CWL$ is not idempotent.

\label{CWL}
\end{proposition}

\begin{proof} Consider a non-trivial map $\Sigma^{n-2} B\Z /2\ra
K(H,n-1)$, which exists because $H$ is a finite 2-group, and
observe that it remains essential when composing with the covering
map $K(H,n-1)\ra X$ (otherwise it would lift non-trivially to the
discrete space $\Omega BG\simeq G$). On the other hand, as $BG$ is
aspherical and $n>2$, there are no non-trivial maps from
$\Sigma^{n-2} B\Z /2$ to $BG$. Hence, $\textrm{map}_*(B\Z /2,BG)$
is not equivalent to $\textrm{map}_*(B\Z /2,X)$, and so $CWX\ncong
CWBG$. As $X=LCWLBG$ and $BG\simeq LBG$, we are done.
\end{proof}

So our goal now is to check that $LCW$ is not
idempotent either. To undertake this task, we will need to
explicitly compute $CWX$. A theorem of Chach\'olski
\cite[20.10]{Chacholski96} establishes that the
$A$-cellularization of any space $Y$ is the fibre of the $\Sigma
A$-nullification of the cofibre of the map $\bigvee_{[A,Y]_*}A\ra
Y$, where the wedge is extended to all the homotopy classes of
maps $A\ra Y$. In our case then we consider the cofibre sequence
$$
\bigvee_{[B\Z /2,X]_*}B\Z /2\longrightarrow X\longrightarrow C,
$$
and we need to
describe its $\Sigma B\Z /2$-nullification.

To accomplish this we first show that $C$ is simply connected.
Because $G$ is generated by elements of order
two (it is simple), by the Mayer-Vietoris sequence this amounts to
checking that every element of order two in $\pi_1X=G$ is in the
image of the homomorphism induced at the level of fundamental
groups by a certain map $B\Z/ 2\ra X$. So, let $f$ be such a map,
and $x$ be an element of $G$ of order two which is in the image of
$\pi_1f$. As the fibre sequence \ref{looping} is principal and $S^{n+1}\vee M$ is equivalent
to a suspension, the fibration $K(H,n-1)\ra X\ra BG$
which defines $X$ is principal by \cite[Theorem 3.1]{Dwyer09}. Moreover, it is classified by the
cohomology class defined by the composition $BG\ra
B\Gamma^{\wedge}_2\ra K(H,n)$, where the second map is
$(n+1)$-Postnikov section. According to the definition of
$\Gamma$, the composition $B\Z /2\ra BG\ra K(H,n)$ is null-homotopic
for every map $B\Z /2\ra BG$.  This implies in particular that
$f$ lifts to $X$, so $C$ is simply connected.

Now we can describe $C$ more precisely: As $\bigvee B\Z /2$ and
$X$ are rationally trivial, $C$ is too by Mayer-Vietoris sequence
in rational homology. Moreover, as the universal covering
fibration of $X$ is principal, it is preserved by $p$-completion
(\cite[II.5.2]{Bousfield72}), and hence the ``formal completions"
$(\prod_p X^{\wedge}_p)_{\mathbb{Q}}$ and $(\prod_p
C^{\wedge}_p)_{\mathbb{Q}}$ are also trivial. Now, Sullivan's
arithmetic square \cite[VI.8.1]{Bousfield72} implies that the
cofibre splits as a product $\prod C^{\wedge}_q$ over all primes
$q$. Clearly $\textrm{map}_*(\Sigma B\Z /2,C^{\wedge}_q)$ is
weakly trivial for every $q$ odd; moreover, as $\bigvee_{[B\Z
/2,X]_*}B\Z /2$ is $q$-acyclic, the $q$-completion of the cofibre sequence that defines $C$
gives an equivalence $BG^{\wedge}_q\simeq C^{\wedge}_q$.
Then by Chach\'olski's Theorem $CWX$ is equivalent to the
homotopy fibre of $X\ra \mathbf{P}_{\Sigma B\Z
/2}(C^{\wedge}_2)\times \prod_{q\neq 2} BG^{\wedge}_q$, and we
only need to describe $\mathbf{P}_{\Sigma B\Z /2}(C^{\wedge}_2)$
to get a complete characterization of $CWX$.
%
%
%

\begin{proposition}

The space $\mathbf{P}_{\Sigma B\Z /2}(C^{\wedge}_2)$ has the
homotopy type of $B\Gamma^{\wedge}_2$.

\end{proposition}

\begin{proof}


We follow the same ideas as in the proof of \cite[Proposition
5.5]{Flores07}, see also \cite[Section 3]{Flores09} for a comprehensive introduction
to Zabrodsky Lemma in this context. We denote the $\Sigma B\Z /2$-null space
$\mathbf{P}_{\Sigma B\Z /2}(C^{\wedge}_2)$ by $P$.
We need to define
maps between $P$ and $B\Gamma^{\wedge}_2$; and we start by
constructing the map $P\ra B\Gamma^{\wedge}_2$
by considering the diagram
\begin{equation}  \xymatrix{ \bigvee_{[B\Z /2,X]_*}B\Z /2
\ar[r] & X \ar[d]^p \ar[r]^i & C \ar@{-->}[dddl]
\ar[r]^{(-)^{\wedge}_2} & C^{\wedge}_2 \ar@{-->}[dddll] \ar[r]^{\eta} & P \ar@{-->}[dddlll]^g \\
 & BG \ar[d]_{(-)^{\wedge}_2} & & & \\
 & BG^{\wedge}_2 \ar[d]_{B\pi ^{\wedge}_2} & & & \\
 & B\Gamma^{\wedge}_2. & & & } \label{eq:diag1}
\end{equation}
Here $\eta$ denotes the coaugmentation of the $\Sigma B\Z
/2$-nullification, and $B\pi ^{\wedge}_2$ denotes the composition
of the isomorphism $BG^{\wedge}_2\simeq BN^{\wedge}_2$ with the
map induced in $2$-completed classifying spaces by the projection
$N\ra \Gamma$. (Recall that $N$ stands for the normalizer of the
Sylow 2-subgroup $S$ of $G$, which controls fusion in $S$.)

Observe that the composite $B\Z /2\ra X\ra BG$ is trivial for
every map $B\Z /2\ra X$. Therefore, the composite of the vertical maps extends
to $C$ (all the extensions are considered up to homotopy). As
$B\Gamma^{\wedge}_2$ is $\Sigma B\Z /2$-null (see
\cite[9.9]{Miller}) and $2$-complete the corresponding map $C\ra
B\Gamma^{\wedge}_2$ also extends to $C^{\wedge}_2$ and then to
$P$. So we have defined $g$, one of the desired maps.

Again for clarity, we summarize the information of the
second part of the proof in a new diagram, which we prove to be
commutative up to unpointed homotopy:
\begin{equation}  \xymatrix{ \bigvee_{[B\Z /2,X]_*}B\Z /2
\ar[r] & X \ar[d]^p \ar[r]^i & C
\ar[r]^{(-)^{\wedge}_2} & C^{\wedge}_2 \ar[r]^{\eta} & P \\
 & BG \ar@{-->}[ur]_l \ar[d]_{(-)^{\wedge}_2} & & & \\
 & BG^{\wedge}_2 \ar@{-->}[uurr]_{l^{\wedge}_2} \ar[d]_{B\pi ^{\wedge}_2} & & & \\
 & B\Gamma^{\wedge}_2 \ar@{-->}[uuurrr]_f . & & & } \label{eq:diag2}
\end{equation}
To define the second map we consider the covering fibre
sequence:
$$
K(H,n-1)\longrightarrow X \longrightarrow BG.
$$
As $H$ is a 2-torsion finite group and $n>2$, the mapping space
$\textrm{map}_*(K(H,n-1),\Omega P)$ is weakly trivial. But the
composition $K(H,n-1)\ra X\ra C\ra P$ is inessential, so Dwyer's
version of the Zabrodsky Lemma (\cite[Proposition~3.5]{Dwyer96}, see
also \cite[Lemma 2.3]{CCS2})
establishes the existence of an extension (up to unpointed homotopy)
$$
\xymatrix{ X \ar[d] \ar[r] & C \ar[r] & P \\
BG \ar[urr]^l & & }
$$
which, in particular, is null-homotopic if and
only if the upper composition is. As $P$ is 2-complete, the
extension $BG\ra P$ factors through $BG^{\wedge}_2$, which has the
homotopy type of $BN^{\wedge}_2$. Now  for $Z(S) = Cl(S) < S$
as discussed just before Lemma \ref{superperfect}, there is
another fibre sequence $BCl(S)\ra BN\ra B\Gamma$. By the
definition of $\Gamma$ and $P$, the composition
$$
BCl(S)\longrightarrow BN\longrightarrow BN^{\wedge}_2\simeq BG^{\wedge}_2\longrightarrow P
$$
is trivial.  Then as $P$ is $\Sigma B\Z /2$-null, Zabrodsky's Lemma
again shows the existence of an extension from $BN^{\wedge}_2\ra
P$ to $B\Gamma^{\wedge}_2$. This is the other map we needed, and
we denote it by $f$.

Our task now is to establish that $f\circ g$ and $g\simeq f$
are (pointed) homotopic to the respective identities. As $P$ and
$B\Gamma^{\wedge}_2$ are simply connected, it is enough to check
this in the unpointed category. Let us denote by $F$ the composition
of the horizontal maps in the previous diagram, and by $G$ the
composition of the vertical maps.

Consider first the composition $f\circ g:P\ra P$. The universal
property of localization guarantees that $f\circ g\simeq Id_P$ if
and only if $f\circ g\circ\eta\simeq\eta$. Moreover, as $P$ is
2-complete, this is equivalent to establishing that $f\circ
g\circ\eta\circ (-)^{\wedge}_2\simeq\eta\circ (-)^{\wedge}_2$.
Now the Puppe sequence of the cofibration $\bigvee_{[B\Z /2,X]_*}B\Z
/2 \ra X$ implies that the latter is true if and only if the maps
are homotopic when precomposing with $i$. By the definition of
$F$, this is the same as checking that $f\circ g\circ F$ is
homotopic to $F$. According to Diagram~1.1, $g\circ F\simeq G$,
and then Diagram~1.2 implies $f\circ G\simeq F$.

On the other hand, consider $g\circ f:B\Gamma^{\wedge}_2\ra
B\Gamma^{\wedge}_2$. As this space is 2-complete and $\Sigma BZ
/2$-null, Zabrodsky's Lemma again applied to the fibration
$BCl(S)\ra BG\ra B\Gamma$ implies that $g\circ f\simeq
Id_{B\Gamma^{\wedge}_2}$ if and only if $g\circ f\circ
B\pi^{\wedge}_2 \circ (-)^{\wedge}_2$ is homotopic to
$B\pi^{\wedge}_2 \circ (-)^{\wedge}_2$. Applying the Zabrodsky
Lemma again to the universal covering fibration of $X$, we obtain that
the latter is equivalent to $g\circ f\circ B\pi^{\wedge}_2 \circ
(-)^{\wedge}_2\circ p\simeq (-)^{\wedge}_2\circ p$, or, in other
words, $g\circ f\circ G\simeq G$. By Diagram~1.2 now, $f\circ
G\simeq F$, and then by Diagram~1.1, $g\circ F\simeq G$. So we are
done.

\end{proof}

\begin{proposition} The functor $LCW$ is not homotopy idempotent.
\label{LCW}
\end{proposition}

\begin{proof}

The previous proposition gives a fibre sequence
$$
CWX\longrightarrow  X\longrightarrow
B\Gamma^{\wedge}_2 \times \prod_{q\textrm{ }odd}BG^{\wedge}_q
$$
that characterizes the cellularization of $X$.
Recall at this point that $X$ is, by definition, $LCWLBG\simeq
LCWBG$.  We show how this space and the previous
description of $CWX$ give a counterexample for the second part of
Problem 7 of Farjoun in \cite{Mislin 08}.

As the map $X\ra
B\Gamma^{\wedge}_2 \times \prod_{q\textrm{ }odd}BG^{\wedge}_q$
factors through $BG$, its composition with $K(H,n-1)\ra X$ is
homotopically trivial, and thus the long exact homotopy sequence
of the previous fibration defines $\pi_{n-1}CWX$ as an extension
$$\pi_n(B\Gamma^{\wedge}_2 \times \prod_{q\textrm{
}odd}BG^{\wedge}_q)\longrightarrow  \pi_{n-1}CWX\longrightarrow  H.$$
Observe that $L$ does not affect the 2-torsion subgroup of
$\pi_{n-1}CWX$, which is therefore isomorphic to the 2-torsion
subgroup of $\pi_{n-1}LCWX$.  In particular it is \emph{not}
isomorphic to the 2-torsion subgroup of $\pi_{n-1}X=H$ because
$\pi_n(B\Gamma^{\wedge}_2 \times \prod_{q\textrm{
}odd}BG^{\wedge}_q)$, by definition of $n$, is nontrivial. Hence
$LCWBG$ is not homotopy equivalent to $LCWX$, which is by
definition $LCWLCWBG$. So $LCW$ is not homotopy idempotent
either, and the argument is complete.

\end{proof}







\bigskip

\bigskip\noindent
Ram\'on J. Flores\\
Departamento de Estad\'\i stica, Universidad Carlos III
de Madrid, C/ Madrid 126\\
E -- 28903 Getafe\\
e-mail: {\tt rflores@est-econ.uc3m.es}


\begin{thebibliography}{noteFarjoun}

\bibitem[Bou94]{Bousfield75}
A.K. Bousfield.
\newblock The localization of spaces with respect to homology.
\newblock {\em Topology}, 14: 133--150, 1975.

\bibitem[Bou94]{Bousfield94}
A.K. Bousfield.
\newblock Localization and periodicity in unstable homotopy theory.
\newblock {\em J. Amer. Math. Soc.}, 7(4): 831--873, 1994.

\bibitem[Bou97]{Bousfield97}
A.K. Bousfield.
\newblock Homotopical localization of spaces.
\newblock {\em Amer. J. Math.}, 119(6): 1321--1354, 1997.

\bibitem[BK72]{Bousfield72}
A.~K. Bousfield and D.~M. Kan.
\newblock {\em Homotopy limits, completions and localizations}.
\newblock Springer-Verlag, Berlin, 1972.
\newblock Lecture Notes in Mathematics, Vol. 304.

\bibitem[CCS]{CCS2}
N.~Castellana, J.~A. Crespo, and J.~Scherer.
\newblock {P}ostnikov pieces and ${B}{\mathbb {Z}}/p$-homotopy theory.
\newblock {\em Trans. Amer. Math. Soc.}, 359(4):1791--1816, 2007.

\bibitem[Cha96]{Chacholski96}
W.~Chach{\'o}lski.
\newblock On the functors {$CW\sb A$} and {$P\sb A$}.
\newblock {\em Duke Math. J.}, 84(3):599--631, 1996.

\bibitem[DF09]{Dwyer09}
W.G. Dwyer and E.Dror-Farjoun
\newblock Localization and cellularization of principal fibrations.
\newblock, Preprint, available at http://www.nd.edu/~wgd/Dvi/LocalizationCellularizationPrincipalFibrations.pdf

\bibitem[Dwy96]{Dwyer96}
W.~G.~Dwyer.
\newblock The centralizer decomposition of {$BG$}.
\newblock In {\em Algebraic topology: new trends in localization and
  periodicity (Sant Feliu de Gu\'\i xols, 1994)}, volume 136 of {\em Progr.
  Math.}, pages 167--184. Birkh\"auser, Basel, 1996.


\bibitem[Far95]{Farjoun95}
E.~Dror-Farjoun.
\newblock {\em Cellular spaces, null spaces and homotopy localization}, volume
  1622 of {\em Lecture Notes in Mathematics}.
\newblock Springer-Verlag, Berlin, 1996.

\bibitem[FF09]{Flores09}
R.J.~Flores and R. Foote.
\newblock The cellular structure of the classifying spaces of finite groups.
\newblock Preprint, available at http://arxiv.org/pdf/0809.4116.

\bibitem[FS07]{Flores07}
R.J.~Flores and J.~Scherer.
\newblock Cellularization of classifying spaces and fusion properties of finite groups.
\newblock {\em J. Lond. Math. Soc.(2)} 76(1):41--56, 2007.

\bibitem[Gor80]{Gor80}
D.~Gorenstein.
\newblock {\em Finite groups}.
\newblock Chelsea Publishing Co., New York, second edition, 1980.

\bibitem[Lev95]{Levi95}
R.~Levi.
\newblock On finite groups and homotopy theory.
\newblock {\em Mem. Amer. Math. Soc.}, 118(567):xiv+100, 1995.

\bibitem[Mil84]{Miller}
H.~Miller.
\newblock The {S}ullivan conjecture on maps from classifying spaces.
\newblock {\em Ann. of Math. (2)}, 120(1):39--87, 1984.

\bibitem[Mis08]{Mislin 08}
\newblock Guido's book of conjectures.
A gift to Guido Mislin on the occasion of his retirement from ETHZ, June 2006. Collected by Indira Chatterji..
\newblock {\em Enseign. Math.},54(2): vols. 1--3, 2008.

\bibitem[RS01]{Rodriguez01}
J.~L. Rodr{\'{\i}}guez and J.~Scherer.
\newblock Cellular approximations using {M}oore spaces.
\newblock In {\em Cohomological methods in homotopy theory (Bellaterra, 1998)},
  volume 196 of {\em Progr. Math.}, pages 357--374. Birkh\"auser, Basel, 2001.

\end{thebibliography}
\end{document}